\documentclass[a4paper]{article}

\usepackage[letterpaper,top=2cm,bottom=2cm,left=3cm,right=3cm,marginparwidth=20mm]{geometry}

\usepackage[fleqn]{amsmath}
\usepackage{amsthm}
\usepackage{amssymb}
\usepackage{mathtools}
\usepackage{graphicx}
\usepackage[colorlinks=true, allcolors=blue]{hyperref}
\usepackage{xparse}
\usepackage{todonotes}
\usepackage{xspace}
\usepackage{etoolbox}
\usepackage[capitalize,noabbrev]{cleveref}
\usepackage{siunitx}
\usepackage{booktabs}
\usepackage{ifthen}
\usepackage{float}
\usepackage{subcaption}
\usepackage{xltabular}
\usepackage{threeparttable}
\usepackage{tikz}
\usetikzlibrary{shapes.geometric, arrows, positioning}

\DeclareDocumentCommand\orderO{o}{\ensuremath{\mathcal{O}\IfValueTF{#1}{\left(#1\right)}{}}}

\usepackage{lineno}
\usepackage{algorithm}
\usepackage{algpseudocode}
\usepackage{authblk}

\providecommand{\keywords}[1]
{
  \small
  \textbf{\textit{Keywords---}} #1
}

\title{Transformer-Based Flow Shop Scheduling Using MILP-Generated Training Data}
\author[1]{R.Wallrath}
\affil[1]{University of Twente, Faculty of Science and Technology, Sustainable Process Technology, Process Design and Optimization, Drienerlolaan 5, 7522 NB Enschede, The Netherlands}
\newboolean{showresponse}
\setboolean{showresponse}{false}
\ifthenelse{\boolean{showresponse}}{
    \DeclareDocumentCommand\reviewerOneComment{mm}{\todo[inline,color=red!50!white]{\textbf{Reviewer 1:} #1 \textbf{Response:} #2}}
    \DeclareDocumentCommand\reviewerTwoComment{mm}{\todo[inline,color=blue!30!white]{\textbf{Reviewer 2:} #1 \textbf{Response:} #2}}
    
}{
    \DeclareDocumentCommand\reviewerOneComment{mm}{}
    \DeclareDocumentCommand\reviewerTwoComment{mm}{}
    
}

\usepackage[authoryear]{natbib}
\date{}
\begin{document}
\maketitle
\begin{center}
    \small\textit{Preprint --- manuscript currently under peer review}
\end{center}
\begin{abstract}
Advances in machine learning (ML) have created new opportunities to complement traditional operations research (OR) methods. In particular, transformer models can capture complex interactions in token sequences by mapping tokens into a high-dimensional embedding space and propagating contextual information via attention. This makes them a candidate to model non-permutation flow shop scheduling with secondary resources as a next-token prediction task, where tokens represent job–machine–secondary resource tuples. For training, mixed-integer linear programming (MILP)-generated schedules are tokenized and used as next-token prediction data. During inference, partial token sequences (prefixes) are randomly generated and completed by the trained transformer through constrained decoding.
A computational study is conducted on a flow shop with 8 jobs, 4 machines, and 3 secondary resources, where jobs are selected from a fixed pool of 20 jobs that is sampled during training and provides the candidates during prefix completion. The transformer achieves better solution quality (smaller makespans) compared to a genetic algorithm (GA), the NEH heuristic, and random search. It is outperformed only by the MILP model and the iterated greedy (IG) heuristic. The study concludes that transformer models can, to some extent, learn patterns from MILP-optimized non-permutation flow shop schedules and that transformer-based scheduling represents an interesting direction for future research, particularly in settings with a fixed, recurring job set.
\end{abstract}


\section{Introduction}
The continuous development of new solution methods is one of the key success factors of operations research (OR). For the well-known NP-hard optimization problem of flow shop scheduling, methods have evolved from simple heuristics to more sophisticated and compute-intensive methods like genetic algorithms (GA), mixed-integer linear programming (MILP), and constraint programming (CP), enabling better solutions in shorter times. Recently, machine learning (ML) methods have become increasingly competitive alternatives to traditional OR methods for flow shop scheduling.
ML is fundamentally based on pattern recognition in large datasets \citep{bishop2006pattern}. The fact that there are general rules to construct good flowshop schedules, such as Johnson's rule \citep{Johnson54} or the NEH heuristic \citep{nawaz1983}, indicates that good solutions share some structural regularities, which ML methods might be able to recognize and exploit.
Furthermore, two properties motivate the exploration of ML methods for flowshop scheduling.
First, they are based on observations of problems through data and do not require first-principle insights into problems. In industrial practice, this can be a key advantage when first principles are unknown, do not exist, or are too complex to formalize into tractable models. Consequently, if white-box methods like discrete-event simulation (DES), MILP, or CP are used, they often result in abstract or simplified models that possess structural validity \footnote{However, advantages, such as results explainability, extrapolation capability in data-sparse settings, rigorous validation, and the possibility to enforce constraints, result from structural validity.} but lack the depth of real-world constraints or the inherent uncertainty of model parameters. ML models, however, are trained directly on raw, uncurated data that reflect the full complexity of real-world problems.
Second, ML models can leverage an efficiency principle to represent complex data. As pointed out in the Johnson-Lindenstrauss lemma \citep{johnson1984extensions}, high-dimensional data can be mapped into a lower-dimensional space while approximately preserving the pairwise distances between datapoints.
This principle enables ML models to map complex datasets into very compact representations without substantial loss of multidimensional structure, which raises the question of whether ML models can also learn to approximate NP-hard optimization problems like flowshop scheduling.

Motivated by advances in large language models (LLMs), this work investigates whether transformer models \citep{vaswani2017attention} can learn the language of optimal flow shop schedules. In a flow shop, operations can be interpreted as tokens, and a schedule as a structured sequence of tokens, analogous to a sentence. LLMs can restore syntactic and semantic coherence when a single word is modified, even if this requires rearranging the whole sentence. Similarly, in flow shops, changing a single job can render the entire schedule suboptimal due to precedence, machine, and secondary-resource interactions and necessitate global reordering. By using attention mechanisms \citep{bahdanau2014neural}, transformer networks have demonstrated strong capability in modeling long-range and high-order interactions in discrete sequences. This makes them a candidate for learning the combinatorial optimization problem of flow shop scheduling.

To date, research on learning-based scheduling has mainly been focused on reinforcement learning (RL) featuring different neural network architectures \citep{wang2023flexible,Gebreyesus2024SpatioTemporal,zhang2020learning}, enhancements of the agent training loop like imitation learning or self-labeling \citep{pan2023knowledge,li2024learning,corsini2024self}, and clever reward or action space design \citep{wang2021dynamic,shahrabi2017reinforcement}. In contrast, this work uses supervised learning \citep{kotary2022fast, ingimundardottir2018discovering}. A decoder-only transformer is trained in next-token prediction on MILP-generated schedules of a non-permutation flow shop. The flow shop is studied on a fixed, recurring job pool and under
limited secondary resources, which reflects individual problem instances, like weekly production planning in a plant over an established product portfolio. The goal of this work is to show at a limited problem scale that a generic sequence model such as the transformer can learn and reproduce patterns of rigorously optimized non-permutation flow shop schedules. 

Next, the paper reviews related work. Thereafter, transformer-based flow shop scheduling is introduced, which is followed by a computational study and a concluding discussion.

\section{Related Work}
\label{bdes2_lit}
A comprehensive review of non-permutation flow shop scheduling solved by traditional OR methods is given in \citep{rossit2018non}. Beyond these methods, ML-based techniques have also been developed and investigated.
In early research, flow shop scheduling problems were mapped onto neural network architectures, following the idea of the Hopfield network \cite{hopfield1985neural}, which suggests that the objective function and problem constraints are encoded in an energy function, and the network dynamics evolve toward stable equilibrium states corresponding to local optima of the scheduling problem \citep{foo1988stochastic,Willems1994Neural}.
However, these early approaches were limited by computational resources and did not use data-driven learning. With advances in ML algorithms and hardware, this situation has changed. Currently, a predominant approach is to use RL agents that interact with scheduling environments. Various neural architectures and training strategies have been proposed, including attention-based actor–critic models, graph neural networks, imitation learning, self-labeling, and hybrid approaches that combine supervised signals with RL.

Several studies focus on pure RL-based training in which the agent learns optimal scheduling decisions only through interaction with the environment.
In \citep{wang2021dynamic}, dynamic adaptive job shop scheduling with unexpected conditions is learned by a deep reinforcement learning (DRL) agent. The resulting makespans are similar to those of heuristic methods (SPT, LPT, FIFO) and a genetic algorithm. A proximal policy optimization (PPO) agent shows stable reward and makespan evolution after 6000 training episodes and is tested on instances up to 10 jobs and 10 machines. To evaluate the actions of the agent immediately, machine utilization is used as the reward function, and the makespan is added only after all jobs are scheduled. The study exemplifies the advantage of ML methods like DRL to learn optimized scheduling strategies in a black-box fashion by interacting with the environment, which contains machine breakdowns and other factors that render generic heuristics less effective. However, it is not discussed how adjusted heuristics or rigorous methods like MILP would perform in comparison to DRL. The authors also point out the inherent challenge of a growing state and action space with increasing problem sizes.
In \citep{zhang2020learning}, job shop scheduling is formulated as a sequential
dispatching problem on a disjunctive graph. A graph neural network policy selects the
next operation to dispatch and is trained end-to-end via PPO with a reward derived from
the change in a makespan lower bound. The learned dispatching rules are size-agnostic
and outperform classical priority dispatching rules.
In \citep{cho2022minimize}, a pointer network is used to minimize the makespan of a permutation flowshop problem of a block assembly process of shipbuilding. The pointer network is trained through a modified RL approach and compared to the SPT, LPT, and NEH heuristics as well as a GA and simulated annealing. The results show that the pointer network is the second-best method, and slightly inferior to the NEH heuristic. In addition, the comparison does not include a rigorous method like MILP.
In \citep{wang2023flexible}, flexible job-shop scheduling is addressed by an actor-critic DRL approach that features a dual attention network (DAN) which has operation and machine message attention blocks to capture relationships between operations and machines. The reward is defined as an estimated lower bound on the makespan to provide feedback for each operation scheduled, and the policy network is trained through PPO. The authors report superior performance compared to four common dispatching rules and comparable performance compared to CP on instances up to 20 jobs and 10 machines. However, a time limit of 30 minutes is enforced on the CP algorithm, and the DRL results are compared to suboptimal solutions in most cases. The training time of the DRL agent is not discussed. Furthermore, on unseen, out-of-the-distribution instances with 10 to 15 jobs, the DRL results are inferior to CP and a GA.
More recently, \citep{xiao2026resched} propose a transformer-based architecture with simplified state representations for flexible job shop scheduling. The policy is, however, again trained through RL and environment interaction rather than through static supervised learning on optimized schedules. Furthermore, secondary resource constraints are not considered.

Beyond pure RL approaches, some studies integrate additional sources of guidance into the training procedure.
In \citep{pan2023knowledge}, the policy network of a DRL agent is trained for flow shop scheduling in a knowledge-guided optimization framework. In this framework, the policy network is trained both through RL and supervised learning. For the latter, poor results of the policy network are improved using a local search operator, which generates labelled data for the supervised learning loop. In the policy network, processing time vectors are embedded and propagated using recurrent, convolutional, and attention layers. Compared to the NEH heuristic and MILP, the approach shows competitive performance for instances with 30 to 100 jobs and five machines. However, the parameters of the MILP solver (Gurobi) are not mentioned, and it is well-known that modifying the default parameters can help with large combinatorial problems. Furthermore, it is not known how good the solution quality in terms of the optimality gap is as the solver is terminated after 3600 seconds. The authors find that the approach is not competitive to a GA.
In \cite{li2024learning}, graph-based imitation learning is proposed for permutation
flow shop scheduling. Jobs are represented as nodes of a graph whose edges encode
processing-time differences. Jobs are
scheduled one after another by an attention decoder trained with a supervised loss to
imitate schedules generated by the NEH heuristic. Hence, the learned policy inherits the solution quality of NEH
rather than that of MILP. In addition, the architecture is a customized graph neural network that requires the definition of nodes, edges, and their
features, while this work uses a standard transformer that operates
directly on operation tokens and encodes no problem-specific structure beyond the
tokenization. While the graph-based imitation approach is promising for mimicking heuristics, it has not been explored how it performs on MILP-generated schedules, which could lead to more complex and difficult-to-learn solution patterns. 
In \citep{corsini2024self}, a self-labeling strategy is proposed for job shop scheduling
problems. A pointer network constructs schedules, and for every
training instance, 256 parallel rollouts are sampled from the currently trained network
by drawing each decision from its output probabilities. The rollout with the minimum
makespan is then used as a pseudo-label to train the network. Hence, the solution quality is inherently bounded by the
self-labeling mechanism, i.e., by the best schedules the model's own search discovers.
This work instead proposes an external optimizer to teach optimality, where the
supervision signal comes from an MILP model.
In \citep{pan2021deep}, a policy network is proposed that embeds processing time vectors using recurrent, fully-connected, and convolutional layers, and uses an attention mechanism for decoding. In a second step, the solutions of the network are improved heuristically using NEH. Compared to longest-processing-time ordering of jobs, random ordering, and variants of the policy network combined with different heuristic improvements, the network yields smaller relative makespan gaps while the computational effort scales similarly to that of heuristics. The authors conclude that satisfactory performance comes from DRL-based generation of feasible schedules that the NEH heuristic improves. However, they also suggest investigating end-to-end algorithms in future work.
Other RL approaches demonstrate that agents can learn to select dispatch rules \citep{ren2021solving,yang2021intelligent,wang2023solving} and operational parameters of heuristics \citep{shahrabi2017reinforcement}, avoiding the need to learn low-level scheduling decisions. While this is an effective approach that exposes agents to reduced combinatorial complexity, the performance of the agents is always tied to the underlying heuristics.
In \citep{ingimundardottir2018discovering}, dispatching rules for the job shop problem are learned by imitation, where a linear preference model is trained to reproduce the
dispatching decisions of a MILP-based expert policy that always dispatches the optimal operation.
While the training labels stem from rigorous optimization, the learned model is a linear priority rule over hand-crafted features rather than a generative sequence model, and the approach targets permutation job shop scheduling without secondary resources.
In \citep{kotary2022fast}, a deep neural network is trained to directly
regress task start times for job shop instances, using training labels from a CP solver under a 30-minute time limit. A post-processing step recovers a feasible schedule from otherwise infeasible predictions. The resulting model outperforms five priority dispatching rules (SPT, LWR, MWR, LOR, MOR) and it takes the CP solver longer to match the solution quality.

Overall, learning-based scheduling research is dominated by RL, imitation learning, and
self-labeling, which rely on reward designs, feature-based encodings, and graph
architectures. Supervised learning approaches that train models directly on rigorously optimized schedules have received comparatively little attention. Herein, this work proposes a decoder-only transformer which is trained with a next-token prediction loss on MILP-generated schedules of a non-permutation flow shop with secondary resource constraints.
\section{Transformer-based flow shop scheduling}
\subsection{Problem description and training data generation}
\label{sec:problem}

\begin{figure}[H]
\centering
\includegraphics[width=\textwidth]{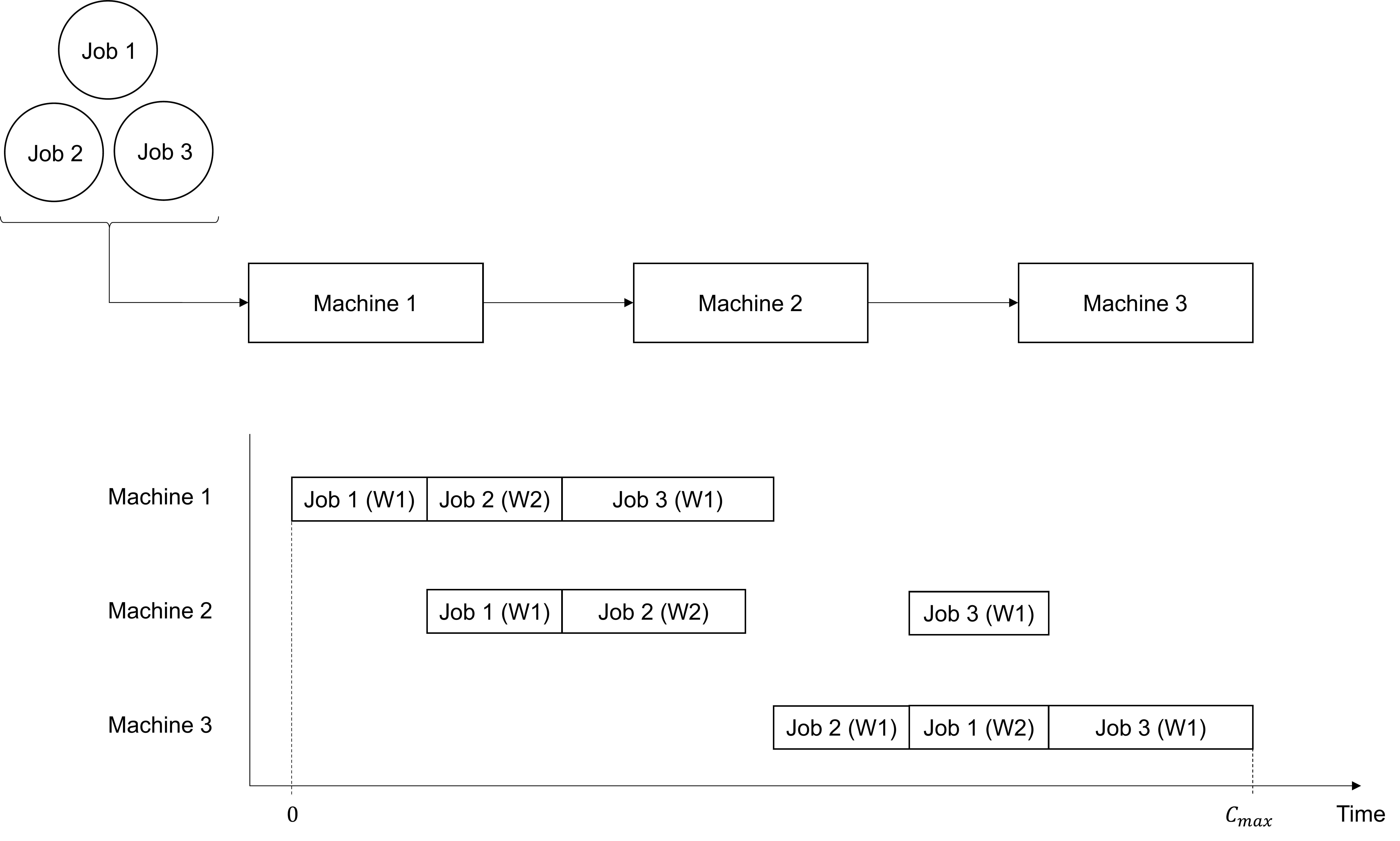}
\caption{Example non-permutation flow shop schedule with $N=3$ jobs, $M=3$ machines, $\mathcal{W}=\{\text{W1}, \text{W2}\}$ with $\text{cap}_1 = \text{cap}_2 = 1$, and job-worker group eligibilities $W_1=W_2=\{\text{W1},\text{W2}\}, W_3 = \{\text{W1}\}$, which causes job 3 on machine 2 to wait for W1.}
\label{fig:flowshop_overview}
\end{figure}

A non-permutation flow shop scheduling problem with $N$ jobs and $M$ machines is considered. All jobs $j \in \{1,...,N\}$ must visit all machines $k \in \{1,...,M\}$ in the same order, while the processing sequence of jobs may differ between machines. Each job has a processing time $p_{jk}$ per machine, and each machine can process at most one job at a time. The objective is to minimize the makespan $C_{\max}$. In addition to machines, jobs require a secondary resource. A global set of worker groups $\mathcal{W}$ with limited capacities $\text{cap}_w$, $w \in \mathcal{W}$, is given. Each job $j$ is compatible only with a subset $W_j \subseteq \mathcal{W}$ of worker groups. When a job starts on machine $k$, it seizes one unit of capacity of an eligible worker group $w \in W_j$ for its entire processing time. To generate a flow shop instance, $N$ jobs are selected from a fixed candidate pool of $N_{\text{pool}}$ jobs with known processing times, eligibilities, and worker group capacities. An example flow shop with $N=3$ jobs, $M=3$ machines, $\mathcal{W}=\{\text{W1}, \text{W2}\}$ worker groups each with $\text{cap}_1 = \text{cap}_2 = 1$, and job-worker group eligibilities $W_1=W_2=\{\text{W1},\text{W2}\}, W_3 = \{\text{W1}\}$ is shown in Figure \ref{fig:flowshop_overview}. In the example, job 3 on machine 2 has to wait for worker group W1 to become available.
\noindent To be used as training data, each flow shop schedule is converted into a sequence of operation tokens $(x_1,\dots,x_{N\cdot M})$, where each token $x_i$ corresponds to a tuple $(j, k, w)$, which encodes the job index $j$, the machine index $k$, and the assigned worker group $w$. Since every one of the $N$ jobs visits all $M$ machines, a complete schedule corresponds to exactly $T = N \cdot M$ operation tokens. The order of the tokens is defined by the operation start times in the schedule. Ties are broken by machine index. A global vocabulary $\mathcal{V}$ of all distinct feasible tuples $(j,k,w)$ over the candidate pool is constructed, and each tuple is assigned to a token index.
To generate the schedules, the MILP model of the flow shop scheduling problem (see Appendix \ref{app:milp}) is solved. First, a job pool of size $N_{\text{pool}}$ with randomized processing times $p_{jk}$, a random job–worker group eligibility matrix, and worker group capacities is created. Second, a subset of $N$ jobs is sampled from the job pool and is solved together with the corresponding eligibility and capacity parameters. This step is repeated $K$ times to obtain a rich training data set containing schedules for many different subsets of the job pool. Each schedule $(x_1,\dots,x_{T})$ is decomposed into multiple training records. For every position $2 \leq t \leq T-1$ in the schedule, the prefix $(x_1,\dots,x_{t-1})$ serves as input and $x_t$ as target, yielding prefix–next-token pairs.
The resulting token sequences are padded to equal length $T$, and attention masks prevent attending to padding positions. The full dataset is divided into a training (90\%) and a validation set (10\%).

\subsection{Transformer-based schedule completion}
\label{sec:transformer_completion}

The transformer model treats the flow shop scheduling problem as a next-token prediction task, that is, given a prefix $(x_1,\dots,x_{t-1})$ of operation tokens,
the transformer learns the conditional distribution
$P(x_{t} \mid x_1,\dots,x_{t-1})$ over all vocabulary tokens $\mathcal{V}$.
At inference time, the model completes a given prefix of length $t_p$ autoregressively. The model predicts the next token, which is appended to the sequence, and the process is repeated until the target length $T$ is reached.
\paragraph{Constrained decoding.}
To guarantee feasibility and full completion, each predicted token is checked against a dynamically updated schedule state, and infeasible tokens are masked out. The state tracks, for every job $j$, the index of the next required machine and the set of jobs already opened. A candidate token $(j,k,w)$ is feasible only if (i) $k$ equals the next required machine of job $j$ in the flow shop, (ii) $w \in W_j$ is an eligible worker group, (iii) job $j$ is not yet finished, and (iv) opening job $j$ does not exceed the maximum of $N$ distinct jobs in the schedule. Rules (i) and (iii) prevent the same operation from appearing twice, and rule (iv) limits the number of distinct jobs to $N$. Since the feasibility check is applied at every decoding step, every emitted sequence is, by construction, a valid schedule.

\begin{figure}[H]
\centering
\includegraphics[width=\textwidth]{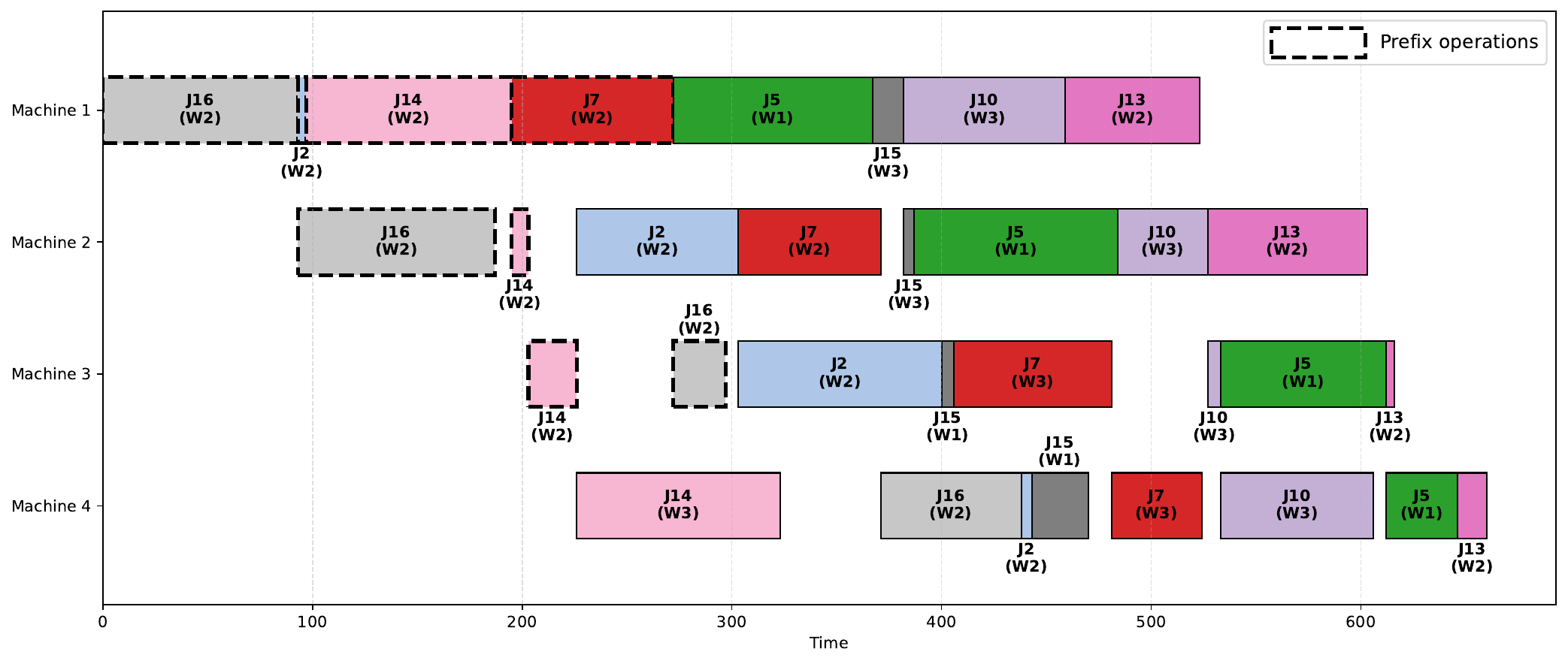}
\caption{Transformer-completed flow shop schedule with prefix of length 8 highlighted.}
\label{Fig:example_gantt}
\end{figure}

\paragraph{Sampling-based search.}
The transformer is evaluated by sampling-based search. For each prefix, $K_s$ completions are generated. Sample $0$ is the deterministic greedy completion. Samples $1,\dots,K_s-1$ are drawn using two control parameters. Before the softmax computation, token logits are divided by a temperature parameter $T$. Smaller values of $T$ amplify the differences between logits and concentrate probability mass on the most likely tokens. A nucleus parameter $p$ restricts sampling to the smallest set of tokens whose cumulative probability reaches $p$, thereby excluding unlikely tokens. However, different settings for $(T,p)$ showed that $(1,1)$, i.e., using the unmodified output distribution of the transformer, yielded the best results averaged across all prefix lengths. Every sequence is then simulated with the ground-truth DES model, and the sequence with the smallest makespan is reported.
An example of a transformer-completed schedule with $t_p = 8$ is shown in Figure \ref{Fig:example_gantt}.

\section{Computational study}
\label{sec:study}
The transformer model is trained on a flow shop with $N=8$ jobs, $M=4$ machines, and a pool of $N_{pool}=20$ jobs, from which $K=40000$ schedules are generated by sampling the job pool. Processing times are integers drawn uniformly from $[1,100]$. There are $\mathcal{W}=\{\text{W1}, \text{W2}, \text{W3}\}$ worker groups with capacities $(2,2,3)$, and the job-worker group eligibility matrix is determined randomly. Each job is eligible for one, two, or all three worker groups, resulting in 38 eligible (job, worker group) pairs and a token vocabulary of $|\mathcal{V}| = 152$ tuples. The training schedules are generated by an MILP model of the flowshop, which includes the secondary resource constraints and is solved with a time limit of 30 minutes using Gurobi 13.0.0, with \texttt{MIPFocus=1} and \texttt{NoRelHeurTime=300}, on 2 threads @ 2.00\,GHz. The transformer (see Appendix \ref{app:architecture}) is trained for 300 epochs with batch size 32 using AdamW (learning rate $10^{-4}$) to minimize the cross-entropy loss.
\subsection{Evaluation process}
\label{sec:evaluation}
The trained transformer network is benchmarked against the IG heuristic \citep{ruiz2007simple}, the NEH heuristic \citep{nawaz1983}, a GA implemented with PyGAD \citep{gad2024pygad}, an MILP model (see Appendix \ref{app:milp}), and random search on their ability to complete partial schedules (prefixes) of lengths $t_p \in {1,\dots,16}$ while minimizing the makespan. For each prefix length, 50 random prefix instances are generated and completed by all methods. The IG and NEH heuristics are chosen as problem-specific methods. The GA is selected as a well-known metaheuristic for this problem \citep{rossit2018non}. MILP is adopted as a rigorous reference. The random search serves as a performance baseline.

\begin{figure}[H]
    \centering
\usetikzlibrary{positioning,calc}
\begin{tikzpicture}[
node distance=0.8cm and 2.5cm,
process/.style={
    trapezium,
    trapezium left angle=70,
    trapezium right angle=110,
    draw,
    text centered
},
data/.style={
    rectangle,
    rounded corners,
    draw,
    minimum width=5.5cm,
    minimum height=1cm,
    text centered
},
arrow/.style={->,thick}
]
\node[process] (prefixgen) {Random prefix generation};
\node[data, below=of prefixgen] (prefix) {Prefix};
\node[data,
      right=1.0cm of prefix,
      minimum width=3.8cm] (jobpool)
{Job pool};
\coordinate (methodscenter) at ($(prefix.south)+(0,-2.0cm)$);
\node[process] (transformer) at ($(methodscenter)+(-5.75cm,0)$) {Transformer};
\node[process] (ig)          at ($(methodscenter)+(-3.45cm,0)$) {IG};
\node[process] (neh)         at ($(methodscenter)+(-1.15cm,0)$) {NEH};
\node[process] (ga)          at ($(methodscenter)+( 1.15cm,0)$) {GA};
\node[process] (milp)        at ($(methodscenter)+( 3.45cm,0)$) {MILP};
\node[process] (rand)        at ($(methodscenter)+( 5.75cm,0)$) {Random};
\node[data, below=1.8cm of methodscenter] (merge)
{Token sequences};
\node[process, below=of merge] (des)
{DES};
\node[data, below=of des] (mk)
{Schedules with makespan};
\draw[arrow] (prefixgen) -- (prefix);
\draw[arrow] (merge) -- (des);
\draw[arrow] (des) -- (mk);
\coordinate (busL) at ($(transformer.north)+(0,0.45cm)$);
\coordinate (busR) at ($(rand.north)+(0,0.45cm)$);
\draw[thick] (busL) -- (busR);
\coordinate (topmid) at (methodscenter |- busL);
\draw[thick] (prefix.south) -- (topmid);
\draw[thick] (jobpool.south) |- ($(prefix.south)+(0,-0.5cm)$);
\foreach \m in {transformer,ig,neh,ga,milp,rand}
{
    \draw[arrow] ($(\m.north)+(0,0.45cm)$) -- (\m.north);
}
\coordinate (bus2L) at ($(transformer.south)+(0,-0.45cm)$);
\coordinate (bus2R) at ($(rand.south)+(0,-0.45cm)$);
\foreach \m in {transformer,ig,neh,ga,milp,rand}
{
    \draw[thick] (\m.south) -- ($(\m.south)+(0,-0.45cm)$);
}
\draw[thick] (bus2L) -- (bus2R);
\coordinate (botmid) at (methodscenter |- bus2L);
\draw[arrow] (botmid) -- (merge.north);
\end{tikzpicture}
\caption{Evaluation process: Every method receives the same prefix and the full job pool to select and schedule the remaining operations. Computation steps are shown as parallelograms, data sets as rectangles.}
\label{Fig:evaluation}
\end{figure}
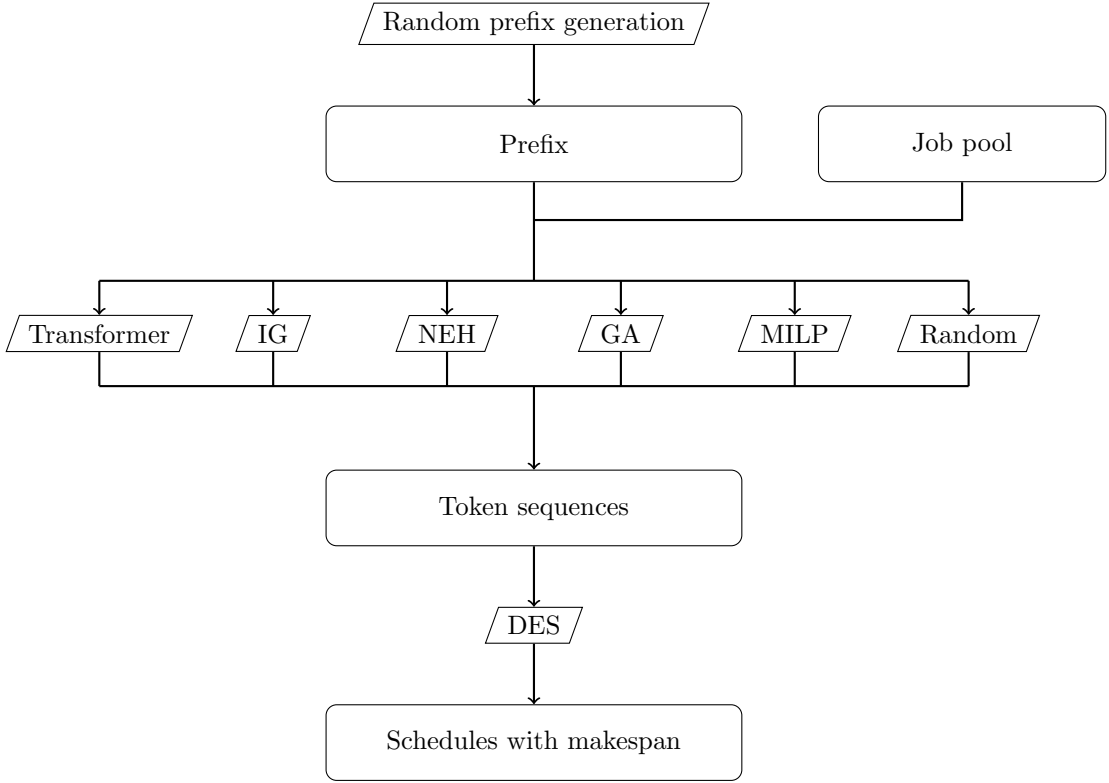

The evaluation process is shown in Figure \ref{Fig:evaluation}. Each method operates in the same decision space, receives the same prefix and job pool, and selects and schedules operations until job completion. In addition, the worker-group assignment of every non-prefix operation is a degree of freedom in every method. Different methods may complete the same prefix with different job sets. However, the requirement to fully complete exactly $N$ jobs makes the comparison on the makespan fair.
Once a method has completed the prefix, the resulting token sequence is rendered into a full schedule using a ground-truth DES model of the flow shop (implemented in SimPy \citep{matloff2008introduction}), and the reported performance measure of every method is the DES makespan of the full schedule. In addition, the internal search of the NEH, IG, and GA is guided by a faster surrogate which is consistent with the DES and models each worker group as $\text{cap}_w$ parallel first-in-first-out slots.

\paragraph{Method implementation.}
Since NEH and IG are originally designed for permutation flow shops, and a general-purpose GA library does not provide scheduling capabilities out of the box, problem-specific implementations are adopted. The transformer, IG, GA, and random search are given the same budget of $B = 1600$ candidate-schedule evaluations per prefix instance to assess their search efficiency and isolate it from implementation and hardware effects. The deterministic NEH heuristic constructs a single schedule and has no budget to scale. The MILP is limited to \SI{600}{\second} solution time. The implementations of all methods are available in a GitHub repository\footnote{\url{https://github.com/roderichwr/Transformer-Based-Flow-Shop-Scheduling}}.

NEH is implemented as an insertion heuristic on the job operations level. All valid positions in the flow shop combined with all eligible worker groups are evaluated, and the best position-group-pair is committed. Operations of different jobs may interleave, which realizes the non-permutation property. The remaining operations of the prefix jobs are inserted first, in descending order of their total remaining processing time. Afterwards, the remaining jobs are selected greedily. Each candidate from the job pool is tentatively inserted at the operation level, and the candidate that minimizes the resulting makespan is committed until $N$ jobs are selected.

IG starts from scratch with a random initial completion of the prefix, i.e., uniform random job selection, operation sequencing, and eligible worker group selection. Each IG iteration randomly chooses $d=2$ jobs and removes all their operations unless they belong to the prefix. Thereafter, the freed operation slots are refilled successively to yield the lowest makespan. If a removed job has operations in the prefix, it is mandatory, and its remaining operations are reinserted at the best (position, worker group) in the completion sequence. The remaining slots are filled at the best (position, worker group) with whichever job operations from the pool minimize the makespan.  
Better solutions are always accepted, and worse ones are accepted with probability 0.1. Since one iteration produces and evaluates exactly one complete candidate schedule, the iteration budget is set to B=1600.

For the GA, a random-key genome with three segments is used. Selection keys determine the jobs from the candidate pool, while prefix jobs are always included. Sequencing keys guarantee feasibility and full completion while allowing job interleaving. Worker group keys assign the workers jointly with selection and sequencing. The fitness is the negative makespan computed by the GA-internal surrogate model of the flow shop. The GA runs with population size $40$ for $39$ generations, so that the initial population plus the offspring evaluated over all generations amount to $B = 1600$ schedules.

Random search draws $B = 1600$ completions per prefix instance. Jobs are randomly sampled, the resulting operations are randomly sequenced, and an eligible worker group is assigned randomly. Thereafter, the candidate schedules are simulated with the DES, and the schedule with the smallest makespan is reported.

The MILP model of Appendix \ref{app:milp} simultaneously performs job selection, sequencing, and worker group assignment. The prefix is enforced by fixing operations on the machine level. A fencing constraint prevents the completion operations from interleaving with the prefix operations. The model contains no warm start and no externally provided bound. The big-M value is computed internally from a schedule of the prefix jobs and the candidates with the smallest processing times. The model is solved with Gurobi (\texttt{TimeLimit} = \SI{600}{\second}, \texttt{MIPFocus=1}, on 8 threads  @ 2.00\,GHz).

\subsection{Problem complexity}
There are many different ways to complete a given prefix, which results in a large decision space.
Consider first the sequencing space for a set of $N$ selected jobs. With a prefix of length \(t_p\), the number of remaining schedules is
\[
S_{\text{seq}} = \left( \prod_{m=1}^{M} (N - t_m)! \right) |W|^{\,M N - t_p},
\quad \text{with} \quad \sum_{m=1}^{M} t_m = t_p,
\]
where \(t_m\) is the number of jobs already scheduled on machine \(m\). For example, if \(t_p = 8\), all 8 jobs are fixed on the first machine (\(t_1=8,\; t_2=t_3=t_4=0\)), and all jobs can use all three worker groups ($|W|=3$), then the sequencing space alone is up to
\(S_{\text{seq}} = (8!)^3 \cdot 3^{24} \approx 1.85 \times 10^{25}\).
In addition, job selection is performed. With $n_{\text{pre}}$ distinct jobs in the prefix, there are $\binom{N_{\text{pool}} - n_{\text{pre}}}{N - n_{\text{pre}}}$ possible job sets, each spanning its own sequencing space, so the total decision space is
\[
S_{\text{total}} = \binom{N_{\text{pool}} - n_{\text{pre}}}{\,N - n_{\text{pre}}\,} \cdot S_{\text{seq}},
\]
which adds a factor ranging from $\binom{12}{0}=1$ up to $\binom{19}{7} = 50388$ to $S_{\text{seq}}$.
While all methods operate in this decision space, training gives the transformer a theoretical advantage over the other methods, as it exposes the transformer to a few schedules before the evaluation. However, this advantage is negligible not only because the entire training solution space contains up to
\(\binom{20}{8} \; (8!)^{4} \; 3^{32} \;\approx\; 6.17 \times 10^{38}\)
feasible flow shop schedules, but also because the space of possible prefixes is very large. Random prefixes force the transformer into states that do not coincide with any training data and thereby prevent full recall. A lower bound on the number of possible prefixes is 
\[
S_{\text{prefix}} \ge \frac{N_{\text{pool}}!}{(N_{\text{pool}} - t_p)!} \cdot |W|^{t_p}
\] which, for example, yields approximately \( 2.03 \times 10^{10}\) possibilities when \(N_{\text{pool}}=20\), \(|W|=3\), and the prefix consists of $t_p=6$ different jobs that are assigned to the first machine. The lower bound is valid because restricting the prefix to contain jobs on the first machine only forms a subset of all possible prefixes of that length. Under this restriction, the calculation simplifies to selecting and sequencing \(t_p\) jobs without replacement from a pool of size \(N_{\text{pool}}\) and assigning each job to one of \(|W|\) worker groups.
Given the size of the training solution and prefix space, the transformer's exposure to $K=40000$ training schedules is relatively small, and learning flow shop scheduling by brute force, i.e., by memorizing the training solution or prefix space, is a negligible effect.
\subsection{Results}
\label{sec:results}
\begin{figure}[H]
\centering
\includegraphics[width=\textwidth]{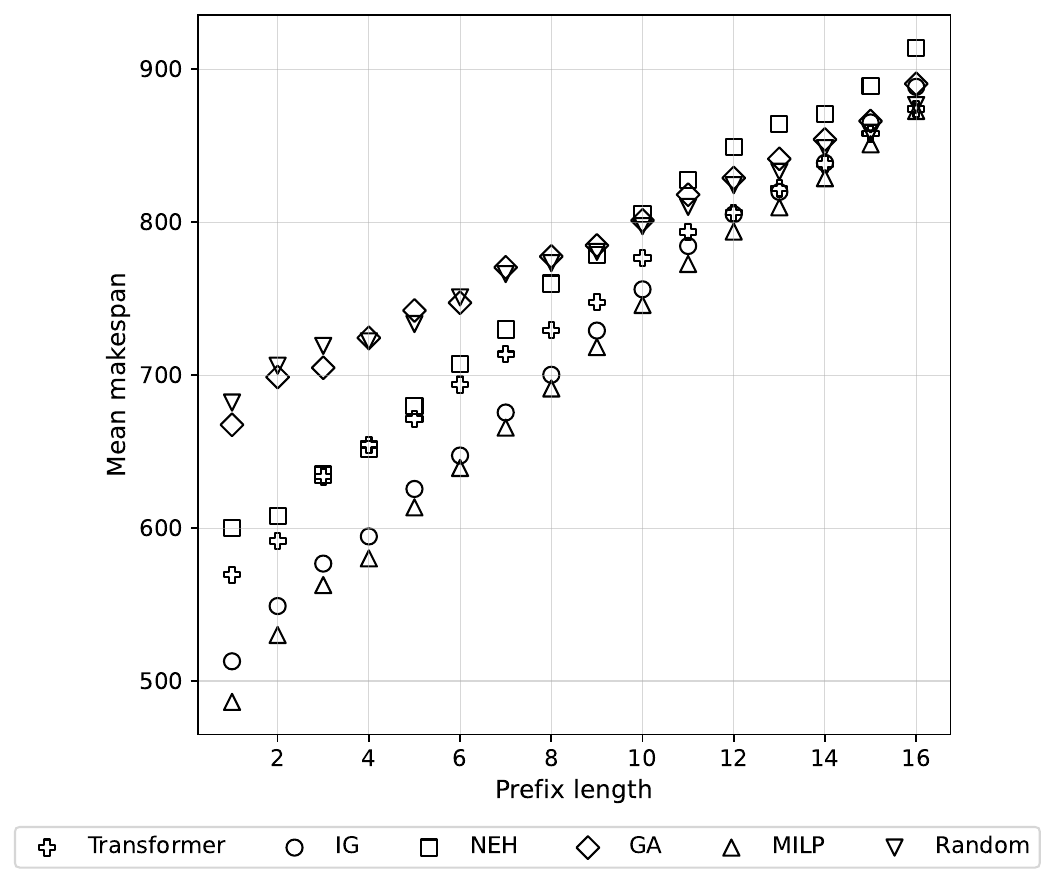}
\caption{Mean DES makespan of 50 randomly generated, partial flow shop schedules (prefixes) completed by the transformer, IG heuristic, NEH heuristic, genetic algorithm, MILP, and random search.}
\label{Fig:completion_rand}
\end{figure}

Figure \ref{Fig:completion_rand} shows the mean DES makespan of 50 random prefixes across prefix lengths $t_p = 1,\dots,16$ completed by the different methods. Detailed results with the distributions are given in Appendix \ref{app:detailed}.
All methods show increasing makespans with increasing prefix length because the portion of randomly fixed operations grows, while the window available for optimized selection and scheduling decisions shrinks. For short prefixes, the difference between the MILP and random search is large,
reaching 40\% at $t_p=1$. Across prefix lengths $t_p = 1,\dots,8$, some of the $N=8$ jobs still have to be selected
from the job pool. Hence, the methods' ability to select makespan-minimizing jobs plays an important role. MILP leads, followed by IG, whose destruction--reconstruction
operator repeatedly re-opens the job set and thus searches selection and sequencing effectively. The transformer and NEH heuristic are, on average, $10.6\%$ and $13\%$ larger than the mean makespans of the MILP. Although the transformer selects jobs only implicitly through its
learned token distribution, it is competitive with the NEH heuristic, which uses an explicit search over candidate jobs. The GA performs similarly to random search, likely due to the limited search budget of 39 generations, which may not be sufficient to improve the schedules. From a
prefix length of $t_p=8$ onwards, the relative performance changes. Since a prefix can contain
at most $N=8$ distinct jobs, the job set is frequently already fixed by the prefix so that scheduling, i.e., the interleaving of the remaining operations and the
assignment of worker groups, becomes the primary tested capability. In this range, the performance gap between the methods narrows from 13.5\% at $t_p=8$ to 4.7\% at $t_p=16$, while the transformer achieves makespans close to those of the MILP model. As the remaining decision space becomes smaller, all methods are able to find near-optimal schedule completions. The random search turns into a competitive alternative and overtakes the
GA, NEH, and at $t_p=15$ the IG, possibly because these methods converge to local optima that
unguided sampling avoids. The random search does not outperform the transformer at any prefix length, although both operate on an
identical budget of 1600 candidate schedules and differ only in the distribution from which
completions are drawn. This indicates that learning has taken place during training and that
the learned distribution concentrates the transformer search on promising regions of the decision space. 
Averaged over all prefix lengths, the transformer finds schedules with mean makespans that are $6\%$ smaller than those of the random search.

\begin{figure}[H]
\centering
\includegraphics[width=\textwidth]{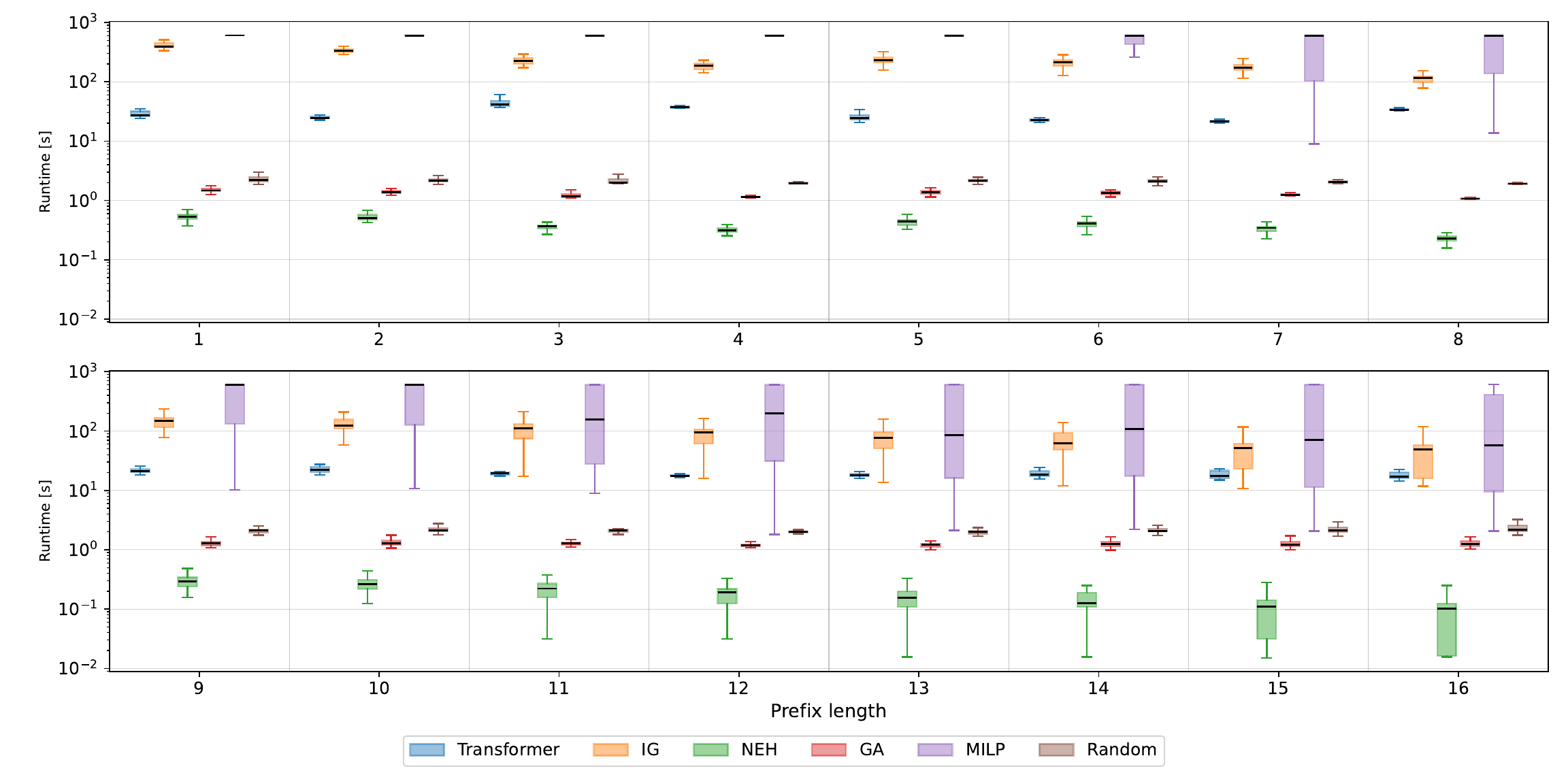}
\caption{Per-instance wall-clock runtime per prefix length for all methods over the 50 evaluated instances (top: $t_p = 1,\dots,8$; bottom: $t_p = 9,\dots,16$).}
\label{Fig:runtime_details}
\end{figure}

\paragraph{Amortized runtime analysis.}
The runtime distributions of the different methods are shown in Figure \ref{Fig:runtime_details}. The MILP and IG take several minutes, while the transformer takes \SI{20} to \SI{40}{\second}, and NEH, GA, and random search compute in a few seconds. However, a complete speed comparison must include the offline costs of the transformer. Generating the $K = 40000$ training schedules with a 30-minute time limit on two threads per MILP-solve requires approximately $40000$ CPU-hours, corresponding to roughly 26 days of wall-clock time on a 64-core machine. Training the transformer added approximately three days. Against MILP, transformer inference saves about \SI{386}{\second} per instance, so the initial investment amortizes after approximately $6500$ scheduling instances. Against IG, which is the closest competitor in solution quality, the saving is about \SI{138}{\second} per instance, and the break-even point is approximately $18000$ scheduling instances. Against NEH, GA, and random search, the transformer does not offer a speed advantage in the current implementation, but its value lies in the quality of the solutions.

\begin{figure}[H]
\centering
\includegraphics[width=\textwidth]{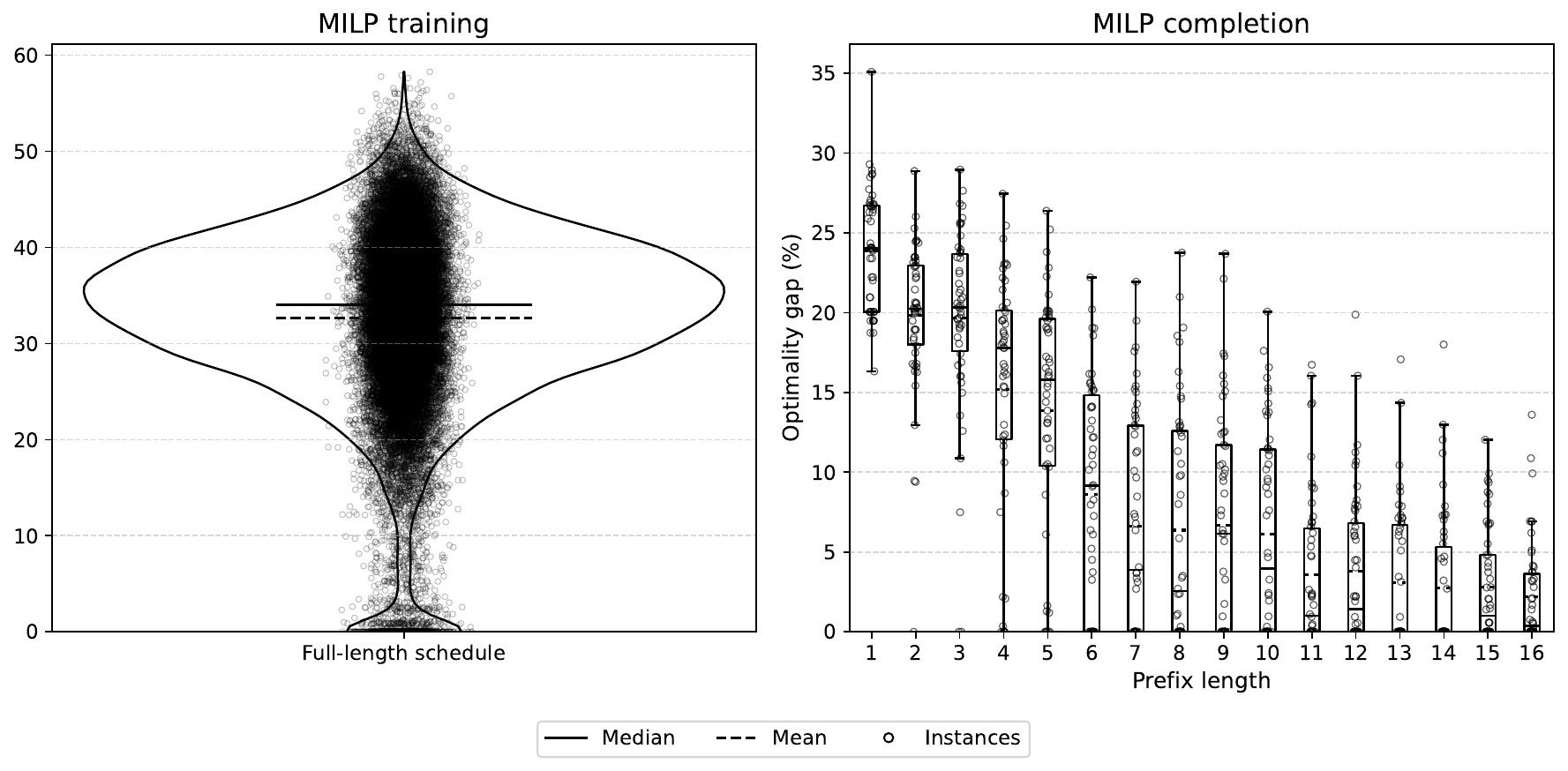}
\caption{Optimality gaps of the MILP model for training data generation with \SI{1800}{\second} (left) and for schedule completion with \SI{600}{\second} time limit (right).}
\label{Fig:MILP_gaps}
\end{figure}

\paragraph{Quality of the training labels.}
The MILP generates full-length training schedules with a 30-minute limit and yields an
optimality gap distribution as shown on the left of Figure \ref{Fig:MILP_gaps}. However, it is known
that MILP models of non-permutation flow shops are difficult to solve to proven
optimality even for small instances \citep{wallrath2025improved}. This is not only due
to the combinatorial size, but also due to big-M reformulations used in the model,
which cause weak dual bounds and slow the convergence. As a result, primal solutions may
already be closer to the optimum than the gaps suggest. The same circumstances apply to
the MILP when used for prefix completion, as shown on the right of Figure
\ref{Fig:MILP_gaps}. Even for the larger gaps at short prefix
lengths, the MILP makespans are the best among all methods (compare
Figure \ref{Fig:completion_rand}).
The gaps decrease as the completion window shrinks, consistent with a dual bound that tightens only slowly against the shrinking combinatorial size of the problem. Independently of the label quality, the supervised learning
concept remains valid. Weaker training schedules would translate into weaker
benchmark performance of the trained model, which is measured directly.

\section{Conclusion}
This study demonstrates that transformer networks can learn non-permutation flow shop scheduling with secondary resources for a fixed job pool as a next-token prediction task through MILP-generated training schedules.
A benchmark is conducted in which random prefixes of varying length are completed by six methods: The transformer, the NEH heuristic, the IG heuristic, a GA, an MILP model, and a random search. Each method has to select the remaining jobs from the same job pool, schedule them on all machines, and assign them to worker groups. The transformer outperforms the GA and the NEH heuristic, and is only outperformed by the IG heuristic and the MILP model. The fact that the transformer outperforms the random search, despite an identical search budget, confirms that learning has taken place.
A runtime analysis shows that search effort is pushed into offline training, while inference is comparatively fast. The offline investment amortizes after $6800$ uses compared to the MILP model, which can be relevant in industrial settings, such as weekly production planning with a recurring product
portfolio.
In addition to solution quality and speed, the transformer approach offers the practical advantage of not requiring a mathematical formulation or a customized (meta)heuristic.
Although training is compute-intensive and requires a large number of sample schedules, this study shows that a relatively simple transformer architecture and training process can already be effective.

The results of this work are preliminary investigations and might not show the full potential of the method. The transformer architecture used is basic and trained under limited computational resources. Larger datasets, hyperparameter tuning, and GPU computations can further improve performance. The model also learns purely from token identities and is therefore tied to the fixed job pool it is trained on. Generalizing to unseen jobs, processing times, or eligibility structures might require feature-based token encodings. Similarly, extending the study beyond the current flow shop size and to an industrial case is a task for future work.

\bibliography{flow-shop-makespan}

\appendix

\section{Transformer architecture}
\label{app:architecture}
Each unique token is embedded via a learnable embedding layer
$\mathbb{R}^{|\mathcal{V}|} \rightarrow \mathbb{R}^{d_{\text{model}}}$ with $d_{\text{model}}=128$.
The transformer model features $L=4$ stacked layers and $H=8$ attention heads per layer. Each layer consists of (i) bidirectional multi-head self-attention with scaled dot-product attention and a padding mask to handle variable-length inputs, (ii) learned relative position bias to encode pairwise distance information, (iii) a position-wise feed-forward network (hidden dimension 256, GELU activation), and (iv) residual connections and layer normalization after each sub-layer.
For an embedded prefix $\mathbf{E} \in \mathbb{R}^{T \times d_{\text{model}}}$, each layer applies self-attention over all prefix positions, enabling the model to capture global dependencies induced by precedence constraints and shared worker resources.
After the final layer, the hidden representation of the last token $\mathbf{h}_t$ is projected via a linear layer to tokens over $\mathcal{V}$.
A softmax produces the probability distribution for the next token.

\section{MILP model}
\label{app:milp}
The MILP model performs job selection, scheduling, and worker group assignment. It is used both as the method in the benchmark and, in an equivalent time-indexed model variant, for training data generation. Let $\mathcal{P} = \{1,\dots,N_{\text{pool}}\}$ denote the candidate pool, $\mathcal{M}=\{1,\dots,M\}$ the machines, and $\mathcal{W}$ the worker groups with capacities $\text{cap}_w$. Each operation is a pair $(j,k)$ with processing time $p_{jk} > 0$. Let $\mathcal{O}_w = \{(j,k) : j \in \mathcal{P},\, k \in \mathcal{M},\, w \in W_j\}$ be the operations eligible for group $w$, and let $H$ be an internally computed scheduling horizon which servers as a valid upper bound obtained from a sequential schedule of the prefix jobs and the smallest-workload candidates.

\begin{align*}
s_j &\in \{0,1\} && \text{Job $j$ is selected,} \\
C_{jk} &\geq 0 && \text{completion time of operation $(j,k)$,} \\
y_{ijk} &\in \{0,1\} && \text{job $i$ precedes job $j$ on machine $k$,} \\
a_{jkw} &\in \{0,1\} && \text{operation $(j,k)$ is served by group $w \in W_j$,} \\
\delta_{oq} &\in \{0,1\} && \text{operation $o$ finishes before operation $q$ starts,} \\
C_{\max} &\geq 0 && \text{makespan.}
\end{align*}

\begin{align}
\min \quad & C_{\max} \\
\text{s.t.} \quad
& \sum_{j \in \mathcal{P}} s_j = N, \qquad s_j = 1 \;\; \forall j \in \text{prefix jobs}, \label{eq:selection}\\
& C_{j1} \geq p_{j1}\, s_j, \qquad C_{jk} \geq C_{j,k-1} + p_{jk} - H(1 - s_j) \quad \forall j, \; k \geq 2, \label{eq:flow}\\
& C_{jm} \leq H s_j \quad \forall j,k, \label{eq:collapse}\\
& C_{ik} \leq C_{jk} - p_{jk} + H(1 - y_{ijk}) \quad \forall i \neq j,\, k, \label{eq:disjunctive}\\
& y_{ijk} + y_{jik} \geq s_i + s_j - 1 \quad \forall i < j,\, k, \label{eq:eitheror}\\
& \sum_{w \in W_j} a_{jkw} = s_j \quad \forall j, k, \label{eq:assign}\\
& C_o \leq C_q - p_q + H(1 - \delta_{oq}) \quad \forall o \neq q, \label{eq:deltalink}\\
& \sum_{\substack{o, q \in Q \\ o \neq q}} \big(\delta_{oq} + \delta_{qo}\big) \;\geq\; \sum_{(j,k) \in Q} a_{jkw} - \text{cap}_w
  \quad \forall w \in \mathcal{W},\; \forall Q \in \mathcal{Q}_w, \label{eq:capacity}\\
& C_{\max} \geq C_{jM} \quad \forall j, \label{eq:makespan}\\
& C_{\max} \geq \sum_{j \in \mathcal{P}} p_{jk}\, s_j \quad \forall k, \qquad
  C_{\max} \geq \frac{1}{\text{cap}_w} \sum_{(j,k) \in \mathcal{O}_w} p_{jk}\, a_{jkw} \quad \forall w. \label{eq:cuts}
\end{align}
Constraint \eqref{eq:selection} selects exactly $N$ jobs and fixes the prefix jobs. Furthermore, the prefix is enforced by fixing the corresponding $y$ variables and worker assignments $a$. A fencing constraint requires all operations of non-prefix selected jobs on a machine to complete after the prefix block on that machine. Constraint \eqref{eq:flow} enforces the flow shop machine order for selected jobs, \eqref{eq:collapse} collapses unselected jobs to $C = 0$, \eqref{eq:disjunctive}–\eqref{eq:eitheror} are the big-M machine disjunctions between selected jobs, and \eqref{eq:assign} assigns every operation of a selected job to exactly one eligible worker group.
Constraint \eqref{eq:capacity} enforces the capacity limit of worker group $w$.
The basic idea is that a group with capacity $\text{cap}_w$ is oversubscribed at
some point in time if and only if at least $\text{cap}_w+1$ of its assigned
operations overlap in time. Consequently, it is sufficient to consider these potentially assigned operation
subsets $Q \in \mathcal{Q}_w$. For each such subset, the constraint requires
that, if all $\text{cap}_w+1$ operations in $Q$ are assigned to worker group
$w$, then at least one pair of them must be sequenced rather than overlap.
The variables $\delta_{oq}$ are linked to the completion times
through \eqref{eq:deltalink}. Specifically, $\delta_{oq}=1$ forces operation
$o$ to finish before operation $q$ starts. Hence, the left-hand side of \eqref{eq:capacity} counts the number of operation pairs in $Q$ that
are explicitly prevented from overlapping.
The right-hand side acts as an activation term. If fewer
than $\text{cap}_w+1$ operations in $Q$ are assigned to worker group $w$, the
right-hand side is non-positive, and the constraint is automatically satisfied.
Only when all operations in $Q$ are assigned to $w$ does the right-hand side
become equal to one, forcing at least one pair of operations to be sequenced.
Constraint \eqref{eq:makespan} defines the makespan. The inequalities in \eqref{eq:cuts} are valid lower bounds
that strengthen the linear relaxation.
The first set of cuts states that the makespan must be at least as large as the
total processing load on each machine. The second set of cuts states that the
makespan must also be at least the total processing load assigned to worker
group $w$, divided by its capacity $\text{cap}_w$.
\section{Detailed makespan results}
\label{app:detailed}
\begin{figure}[H]
\centering
\includegraphics[width=\textwidth]{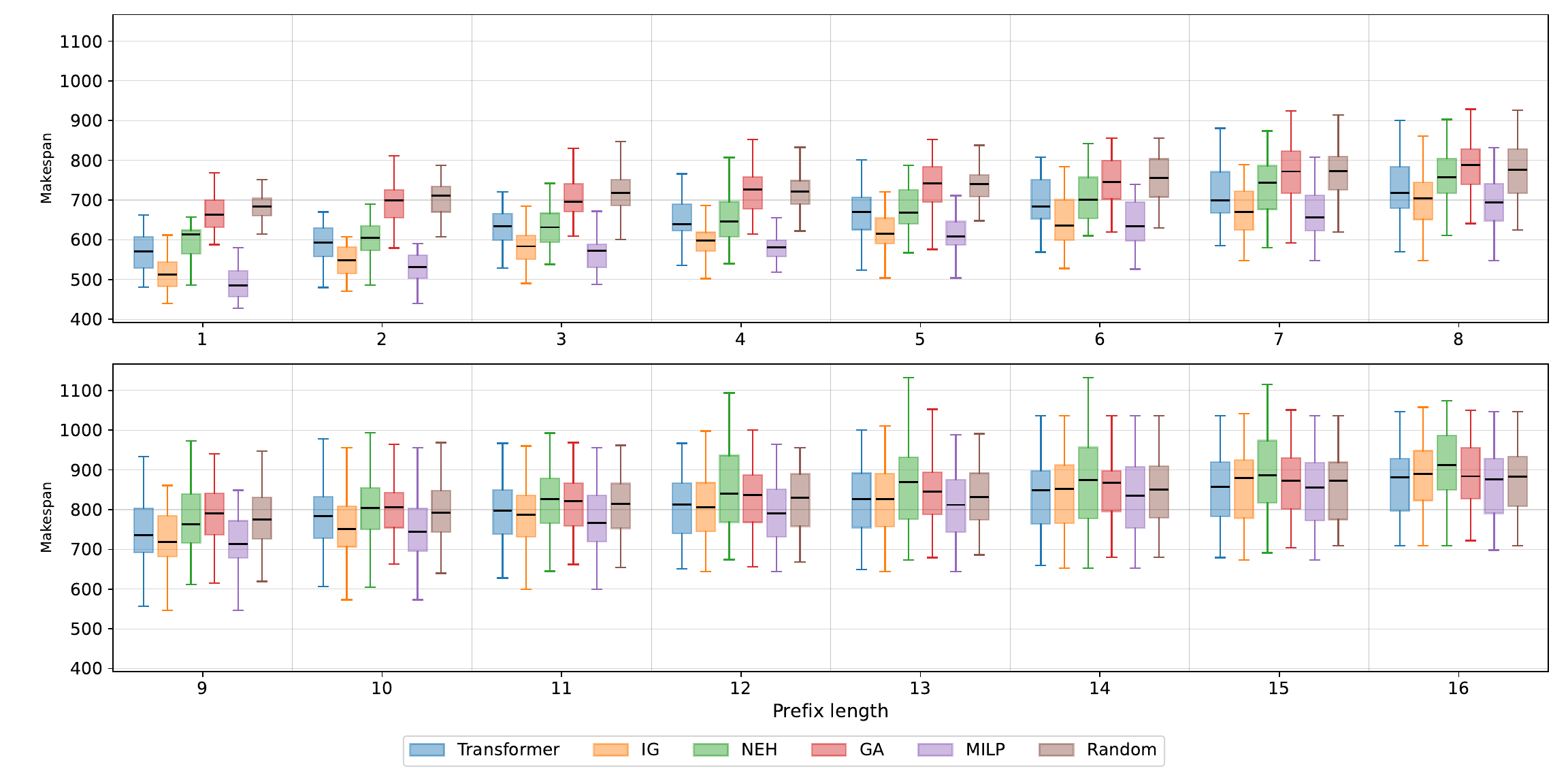}
\caption{DES makespan per prefix length for all methods over the 50 evaluated instances (top: $t_p = 1,\dots,8$; bottom: $t_p = 9,\dots,16$).}
\label{Fig:makespan_details}
\end{figure}

Figure \ref{Fig:makespan_details} shows the distribution of the DES makespans over the 50 instances per prefix length, complementing the means of Figure \ref{Fig:completion_rand}. The spread across instances is considerable, as each seed generates a different random prefix and therefore a different completion problem. However, the relative performance of the methods discussed in Section \ref{sec:results} is consistent across the instance population.
 
\end{document}